\documentclass[twocolumn]{autart}

\usepackage{cite} 		  
\usepackage{hyperref}	  
\usepackage{breakurl}
\hypersetup{colorlinks=true,linkcolor=blue,citecolor=blue,breaklinks={true},allcolors=blue} 
\usepackage[normalem]{ulem}
\usepackage{comment}

\usepackage{graphicx}
\graphicspath{{./FIG/}}
\DeclareGraphicsExtensions{.eps,.png,.jpg,.jpeg,.bmp,.gif,.pdf}

\usepackage{tabularx,booktabs}	
\usepackage{multirow}	        
\usepackage{array}
\newcolumntype{?}{!{\vrule width 2pt}}	
\usepackage{enumitem}           

\usepackage{amsmath}      		
\usepackage{amssymb}	  		
\usepackage{cancel}
\usepackage{bm}		      		
\usepackage{physics}	  		
\usepackage{mathtools}			
\usepackage{dsfont}		  		
\usepackage{soul}				
\usepackage{arydshln}           
\usepackage[version=4]{mhchem}  
\usepackage{multicol}           

\newcommand{\R}{\mathbb{R}}		 
\newcommand{\C}{\mathbb{C}}		 
\newcommand{\transp}{\mathsf{T}}					
\newcommand{\row}{\operatorname{row}}			    
\newcommand{\col}{\operatorname{col}}			    
\newcommand{\diag}{\operatorname{diag}}			    

\newcommand*{\QEDS}{\hfill\ensuremath{\square}}  
\newcommand*{\QEDA}{\hfill\ensuremath{\triangle}}

\newcommand{\rev}[1]{\textcolor{black}{#1}}

\allowdisplaybreaks
\def\arraystretch{1.2}

\usepackage[semicolon]{natbib}
\begin{document}

\setlength{\abovedisplayskip}{7pt}
\setlength{\belowdisplayskip}{7pt}

\begin{frontmatter}

\title{On the Popov-Belevitch-Hautus tests for \\ functional observability and output controllability}

\thanks{Corresponding author: Arthur N. Montanari. Tel. +1 847-467-2552. E-mail: {arthur.montanari@northwestern.edu}}

\vspace{-0.75cm}
\author[1,2]{Arthur N. Montanari}, 
\author[3]{Chao Duan}, 
\author[1,2,4,5]{Adilson E. Motter}, 

\address[1]{Department of Physics and Astronomy, Northwestern University, Evanston, IL 60208, USA} 
\address[2]{Center for Network Dynamics, Northwestern University, Evanston, IL 60208, USA}
\address[3]{School of Electrical Engineering, Xi'an Jiaotong University, Xi'an 710049, China}
\address[4]{Department of Engineering Sciences and Applied Mathematics, Northwestern University, Evanston, IL 60208, USA}
\address[5]{Northwestern Institute on Complex Systems, Northwestern University, Evanston, IL 60208, USA}

\vspace{-0.8cm}
\begin{abstract} \vspace{-0.6cm}

\noindent
Functional observability and output controllability are properties that establish the conditions for the partial estimation and partial control of the system state, respectively. In the special case of full-state observability and controllability, the Popov-Belevitch-Hautus (PBH) tests provide conditions for the properties to hold based on the system eigenspace. Generalizations of the PBH test have been recently proposed for functional observability and output controllability, but thus far have only been proven valid for diagonalizable systems. Here, we rigorously establish the generalized PBH test for functional observability, extending its validity to a broader class of systems using Jordan decomposition. Likewise, we  determine the class of systems under which the generalized PBH test is sufficient and necessary for output controllability. These results have immediate implications for observer and controller design, pole assignment, and optimal placement of sensors and drivers.

\smallskip\noindent
Published in \textit{Automatica} 174:112122 (2025), DOI: \textcolor{blue}{\href{https://doi.org/10.1016/j.automatica.2025.112122}{10.1016/j.automatica.2025.112122}}.

\vspace{-0.5cm}
\end{abstract}

\begin{keyword}  
functional observability; output controllability; Hautus lemma; duality.
\end{keyword} 

\end{frontmatter}

\section{Introduction}

Functional observability is a property of a dynamical system that enables the \textit{partial} reconstruction of the system state from its input and output signals \citep{Fernando2010,Jennings2011}. For linear time-invariant (LTI) systems, this property constitutes a necessary and sufficient condition for the design of functional observers capable of asymptotically estimating a linear \textit{functional} $\bm z(t) = F\bm x(t)\in\R^{r}$ of the system state $\bm x\in\R^n$ \citep{Darouach2000,Fernando2010}. Here, $\bm z$ is referred as the target vector, which typically has a much smaller dimension than the system size (i.e., $r\ll n$). Moving from applications in state estimation to state control, the related notion of output controllability (also known as target controllability) determines the minimal conditions for the existence of a control signal $\bm u(t)$ that steers a target vector $\bm z(t) = F\bm x(t)$ 
to any final state \citep{Bertram1960,Gao2014,Lazar2020}. 

These properties are crucial for solving problems in partial state estimation and partial state control. However, until recently, the development of conditions and methods for the analysis of functional observability \citep{Fernando2010,Jennings2011,Rotella2016,Montanari2022nonlinear} and output controllability \citep{Morse1971,Lazar2020,Danhane2023}
has been pursued independently in the literature.
%
In the case of functional observability, the analysis often relies on rank-based conditions that generalize Kalman's classical observability test \citep{Jennings2011,Rotella2016} or the Popov-Belevitch-Hautus (PBH) test \citep{Moreno2001,Jennings2011}.
This has led to recent advances in algorithms for observer design \citep{Niazi2020,Darouach2020,Darouach2023functdetect,niazi2023clustering}, sensor placement \citep{Montanari2022,Zhang2023b}, and fault/cyberattack detection \citep{emami2015functional,eSousa2022,Venkateswaran2021}. Recent progress in output controllability, on the other hand, has been focused on designing graph-theoretic and cost-effective algorithms for optimal driver placement \citep{Gao2014,Wu2015,Waarde2017,Czeizler2018,Vosughi2019,Li2020,Casadei2020,Li2021,Li2023,zhang2024reachability}. Despite these developments, only recently the duality between functional observability and output controllability was rigorously established \citep{Montanari2023duality,Montanari2023target}, enabling techniques and results derived for one problem to be applied to another (dual) problem. Specifically, functional observability of a system implies the output controllability of the dual system, but the converse only holds under a particular condition.

In the case of full-state estimation/control, there exists a general equivalence between the rank-based tests proposed by Kalman and by Popov, Belevitch, and Hautus. Crucially, such equivalence does not hold for all LTI systems in the case of partial estimation/control. In particular, the generalization of the PBH test for functional observability proposed in \citep{Moreno2001,Jennings2011} has thus far been proved to be valid only for diagonalizable systems \citep{Zhang2023onfunctobsv}; whether such diagonalizability condition can be relaxed is an outstanding question considered here. Moreover, it can be shown that~the recently proposed PBH test for output controllability \citep{Schonlein2023outputctrb} fails for certain nondiagonalizable systems (see Example~\ref{examp.cond3isonlynecessary} below). 
Given the duality between functional observability and output controllability, one might expect that the classes of systems for which the generalized PBH test is valid (or invalid) for each problem would be related. This problem is of paramount relevance to determine the scope of applicability of methods and theoretical results reliant on these tests, including conditions for functional detectability \citep{Darouach2023functdetect,Zhang2023onfunctobsv}, structured systems \citep{Montanari2022,Montanari2023target,Zhang2023onfunctobsv}, and other applications \citep{Niazi2020,eSousa2022,Zhang2023b,niazi2023clustering}.

In this paper, we investigate the conditions for functional observability and output controllability. In particular, we establish a more general class of LTI systems for which the PBH tests for functional observability (Section~\ref{sec.functionalobsv}) and output controllability (Section \ref{sec.outputctrb}) are valid, relaxing the previous diagonalizability condition. We also show that these properties provide \textit{minimal} conditions for the estimation/control of a low-dimensional vector $\bm z(t)$. These properties are inherently related to the observability and controllability Gramians and their ability to reconstruct and steer the target state $\bm z(t)$,~respectively. The duality between functional observability and output controllability allows us to establish an equivalence between the results in Sections \ref{sec.functionalobsv} and \ref{sec.outputctrb}, as shown in Section~\ref{sec.duality}. 
Our identification of a broader class of systems under which the PBH tests are valid has immediate implications for the literature on theory and method development based on these tests.


\section{Preliminaries}
The following notation is adopted. 
The observability matrix of a pair $(M,N)$ is denoted by $\mathcal O(M,N) = [M^\transp \,\, (MN)^\transp \,\, \ldots \,\, (MN^{n-1})^\transp]^\transp$ and the controllability matrix of a pair $(M,N)$ is denoted by $\mathcal C(M,N) = [N \,\, MN \,\, M^2N \,\, \ldots \,\, M^{n-1}N]$.
For compactness, we denote $\mathcal O(C,A)$ and $\mathcal C(A,B)$ simply by $\mathcal O$ and $\mathcal C$, respectively.
%
Let $I_n$ denote an $n\times n$ identity matrix and $0_{m\times n}$ denote an $m\times n$ null matrix (subscripts are omitted when self-evident). 
The block diagonal matrix with submatrices $M$ and $N$ along its diagonal is denoted by $\diag(M,N)$. Likewise, $\diag(\{M_i\}_{i=1}^n)$ denotes a block diagonal matrix formed by submatrices $M_1,\ldots,M_n$. 

The Jordan decomposition of a matrix $A$ is determined by $J = PAP^{-1}= \diag(\{J_k\}_{k=1}^m)$, where each submatrix $J_k = \diag(\{J_{ki}\}_{i=1}^{m_k})\in\R^{n_k \times n_k}$ consists of \textit{all} $m_k$ Jordan blocks $J_{ki}$ associated with the eigenvalue $\lambda_k$ of algebraic multiplicity $n_k$ and geometric multiplicity $m_k$. Let $\bar C = CP^{-1} = [C_1 \,\, \ldots \,\, C_m]$ be partioned as in $J$, where each submatrix $C_k$ is formed by columns corresponding to $J_k$. 
%

\begin{defn}
    For each pair of submatrices $J_k$ and $C_k$ associated with an eigenvalue $\lambda_k$ of multiplicity $n_k>1$, the \textit{lead columns} of $C_k$ are those associated with the first column of each block $J_{ki}$, $i=1,\ldots,m_k$. If $n_k = 1$, $C_k$ does not contain lead columns.
\label{def.jordanleadingcols}
\end{defn}

\begin{lem} \label{lemma.jordanleadingcols}
    If all lead columns of $C_k$ are linearly independent (LI), then
    \begin{equation}
        \rank\begin{bmatrix}
            \lambda_k I_{n_k} - J_k \\ C_k
        \end{bmatrix}
        = n_k
    \label{eq.leadingcolsLI}
    \end{equation}
    \noindent
    and the pair $(C_k,J_k)$ is observable.
\end{lem}

\begin{exmp}
Consider the following system:
\begin{align*}
    J &= \diag(J_1,J_2,J_3) = \left[\begin{array}{cccc|ccc|c}
        \lambda_1 & 1 & 0 & 0 & 0 & 0 & 0 & 0 \\
        0 & \lambda_1 & 1 & 0 & 0 & 0 & 0 & 0 \\
        0 & 0 & \lambda_1 & 0 & 0 & 0 & 0 & 0 \\
        0 & 0 & 0 & \lambda_1 & 0 & 0 & 0 & 0\\
        \hline 
        0 & 0 & 0 & 0 & \lambda_2 & 1 & 0 & 0\\
        0 & 0 & 0 & 0 & 0 & \lambda_2 & 0 & 0\\
        0 & 0 & 0 & 0 & 0 & 0 & \lambda_2 & 0 \\
        \hline 
        0 & 0 & 0 & 0 & 0 & 0 & 0 & \lambda_3
    \end{array}\right],
    \\
    \bar C &= \begin{bmatrix} C_1 & C_2 & C_3 \end{bmatrix} = \left[\begin{array}{cccc|ccc|c} \bm c_{11} & \bm c_{12} & \bm c_{13} & \bm c_{14} & \bm c_{21} & \bm c_{22} & \bm c_{23} & \bm c_{31}\end{array}\right],
\end{align*}
\noindent
where $J$ is in Jordan form. By Definition~\ref{def.jordanleadingcols}, the lead columns of $C_1$ and $C_2$ are $\{\bm c_{11},\bm c_{14}\}$ and $\{\bm c_{21},\bm c_{23}\}$, respectively. The submatrix $C_3$ has no lead columns since $n_k = 1$. Condition \eqref{eq.leadingcolsLI} is satisfied for $k=1$ if the lead columns $\{\bm c_{11},\bm c_{14}\}$ are LI, which implies that the pair $(C_1,J_1)$ is observable.
\QEDA
\end{exmp}

\section{\rev{PBH test for functional observability}}
\label{sec.functionalobsv}

Consider the LTI dynamical system
\begin{align}
    \dot{\bm x} &= A\bm x + B\bm u, 
    \label{eq.dynsys}
    \\
    \bm y &= C\bm x,
    \label{eq.output}
\end{align}
where $\bm x\in\R^n$ is the state vector, $\bm u\in\R^p$ is the input vector, $\bm y\in\R^q$ is the output vector, $A\in\R^{n\times n}$ is the system matrix, $B\in\R^{n\times p}$ is the input matrix, and $C\in\R^{q\times n}$ is the output matrix. 
Let the linear function of the state variables
\begin{equation}
    \bm z = F\bm x
    \label{eq.target}
\end{equation}
define a \textit{target vector} $\bm z \in \R^r$ \textit{sought to be estimated}, where $F\in\R^{r\times n}$ is the functional matrix and $1\leq r\leq n$.

We now define the notion of functional observability and establish the minimal conditions under which a system is functionally observable.

\begin{defn}
\label{def.functobsv}
    \rev{Let $\bm x(t)$ and $\tilde{\bm x}(t)$ be solutions of Eq.~\eqref{eq.dynsys} corresponding to the initial states $\bm x(0)$ and $\tilde{\bm x}(0)$, respectively, under the same control input $\bm u(t)$. The system \eqref{eq.dynsys}--\eqref{eq.target}, or the triple $(C,A;F)$, is functionally observable if, for any pair $\bm x(0)$ and $\tilde{\bm x}(0)$, $C\bm x(t)=C\tilde{\bm x}(t)$ implies $F\bm x(t)=F\tilde{\bm x}(t)$, $\forall t\in[0,t_1]$.}
\end{defn}

\rev{Intuitively, a system is functionally observable if the target state $\bm z(0)=F\bm x(0)$ can be uniquely determined from knowledge of $\bm y(t)$ and $\bm u(t)$ over a finite time $t\in[0,t_1]$.
We assume that $\rank(F)=r$ without loss of generality since, for any $F'$ that is a linear combination of the rows of $F$, it follows $\bm z':=F'\bm x = \bm\gamma^\transp \bm z$ for some $\bm\gamma\in\R^r$. Thus, $\bm z'$ can be directly inferred from $\bm z$, and hence it suffices to verify the functional observability of $(C,A;F)$.}

\begin{thm}
\label{thm.functobsv} {\normalfont\textbf{(Functional observability)}}~Consider the triple $(C,A;F)$ and the equivalent system $(\bar C,J; \bar F)$, where $J=PAP^{-1}$ is in Jordan form, $\bar C = CP^{-1} \rev{=[C_1 \, \ldots \,C_m]}$, and $\bar F=FP^{-1}\rev{=[F_1 \, \ldots \,F_m]}$. The following statements are equivalent:
    \begin{enumerate}
        \item[1)] the triple $(C,A;F)$ is functionally observable;
        \smallskip

        \item[2)] $\rank\begin{bmatrix} \mathcal O \\ F \end{bmatrix} = \rank(\mathcal O)$;
        \smallskip
        
        \item[3)] $\rank 
        \begin{bmatrix}
            \mathcal O \\ \mathcal O(F,A)
        \end{bmatrix}
         = \rank(\mathcal O)$;
        \smallskip
    \end{enumerate}
    \noindent
    Moreover, under the assumption that the lead columns of $F_k$ are LI, for all $k$ such that $F_k\neq 0$, the following statement is also equivalent:
    \begin{enumerate}
        \item[4)] $\rank \begin{bmatrix} \lambda I_n - A \\ C \\ F \end{bmatrix}
        =
        \rank \begin{bmatrix} \lambda I_n - A \\ C \end{bmatrix}$, $\forall\lambda\in\C$.
    \end{enumerate}
\end{thm}

\begin{pf}
%
%
\textit{(1)$\Leftrightarrow$(2)}. \rev{Let $\bm x(t)$ be the solution of system \eqref{eq.dynsys} for an initial state $\bm x(0)$ and an input $\bm u(t)$. Thus,}
\begin{equation}
    \begin{aligned}
    \bm y(t) &= C\bm x(t), \\
    \bm y(t) &= Ce^{At}\bm x(0) + C\int_0^{t} e^{-A\tau}B\bm u(\tau){\rm d}\tau, \\
    Ce^{At}\bm x(0) &= \bm{h}(t,\bm{u},\bm{y}),
    \label{eq.proof.x0y}
    \end{aligned}
\end{equation}
where $\bm{h}(t,\bm{u},\bm{y}) =\bm{y}(t)- \int_0^t e^{-A \tau} \bm{B} \bm{u}(\tau) {\rm d}\tau$ is a functional of  $\bm{u}(t)$ and $\bm{y}(t)$. Multiplying Eq.~\eqref{eq.proof.x0y} on the left by $e^{A^\transp t}C^\transp$ and integrating over $t\in [0,t_1]$ yields
\begin{equation}
\label{eq.proof.multiplybyG}
    W_o(t_1)\bm x(0) = \int_0^{t_1} e^{A^\transp t}C^\transp \bm{h}(t,\bm{u},\bm{y})  \text{d}t,
\end{equation}
where $W_o(t)=\int_0^t e^{A^\transp\tau}C^\transp C e^{A\tau}{\rm d}\tau$ is the observability Gramian. 
\rev{Likewise, given the solution $\tilde{\bm x}(t)$ corresponding to some other initial state $\tilde{\bm x}(0)$ and the same input $\bm u(t)$, we have that}
\begin{equation}
\label{eq.proof.multiplybyG2} \rev{
    W_o(t_1)\tilde{\bm x}(0) = \int_0^{t_1} e^{A^\transp t}C^\transp \bm{h}(t,\bm{u},\tilde{\bm{y}})  \text{d}t,}
\end{equation}
\rev{where $\tilde{\bm y}(t) := C\tilde{\bm x}(t)$. According to Definition \ref{def.functobsv}, a system is functional observable if $\bm y(t) = \tilde{\bm y}(t)$ implies $\bm F\bm x(t)=F\tilde{\bm x}(t)$. It follows that $\bm y(t) = \tilde{\bm y}(t)$ is equivalent to $W_o(t_1)(\bm x(0) - \tilde{\bm x}(0)) = 0$ since the RHS of Eqs. \eqref{eq.proof.multiplybyG} and \eqref{eq.proof.multiplybyG2} are equal. Therefore, if there exists a matrix $G\in\R^{r\times n}$ such that $GW_o(t_1) = F$, it follows that $GW_o(\bm x(0) - \tilde{\bm x}(0)) = F(\bm x(0) - \tilde{\bm x}(0))$, and hence the system is functionally observable.}


The proof now follows from the fact that $GW_o(t_1) = F$ if and only if $\row (F)\subseteq\row(W_o(t_1))$. Since $\operatorname{Im}(W_o)=\operatorname{Im}(\mathcal O)$, it follows that $GW_o(t_1) = F$ holds if and only if $\row(F)\subseteq\row(\mathcal O)$, which is equivalent to statement 2.

\textit{(2)$\Leftrightarrow$(3)}. The equivalence is proven in \citep[Section I]{Rotella2016}.

\textit{(2)$\Rightarrow$(4)}. Let $\rank(\mathcal O) = n_o < n$. 
We apply the following canonical decomposition to system $(C,A;F)$. Let $Q^\transp = [q_1 \,\, \ldots \,\, q_n] \in\R^{n\times n}$ be a unitary matrix such that the first $n_o$ columns $\{q_1,\ldots, q_{n_o}\}$ lie in the $\row(\mathcal O)$ and $\{q_{n_o+1},\ldots,q_n\}$ are arbitrarily chosen such that $Q$ is nonsingular. Applying the similarity transformation $\bar{\bm x} = [\bm x_c^\transp \,\, \bm x_u^\transp]^\transp = Q\bm x$ to system \eqref{eq.dynsys}--\eqref{eq.target} yields \citep{Chi-TsongChen1999}:
\begin{equation}
\begin{aligned}
 \label{eq.decomposedsys}
     \begin{bmatrix}
         \dot{\bm x}_c \\ \dot{\bm x}_u
     \end{bmatrix}
     &=
     \begin{bmatrix}
         A_c & A_{12} \\ 0 & A_{u}
     \end{bmatrix}
     \begin{bmatrix}
         {\bm x}_c \\ {\bm x}_u
     \end{bmatrix}
     +
     \begin{bmatrix}
         B_o \\ B_u
     \end{bmatrix}
     \bm u,
     \\
     \bm y &=
     \begin{bmatrix}
         C_o & 0
     \end{bmatrix}
     \begin{bmatrix}
         {\bm x}_c \\ {\bm x}_u
     \end{bmatrix},
     \\
     \bm z &=
     \begin{bmatrix}
         F_o & F_u
     \end{bmatrix}
     \begin{bmatrix}
         {\bm x}_c \\ {\bm x}_u
     \end{bmatrix},
\end{aligned}
\end{equation}
\noindent
where $\bm x_o\in\R^{n_o}$ and $\bm x_u\in\R^{n-n_o}$ are the observable and unobservable variables, respectively.
Since $FQ^{\transp} = [F_o \,\, F_u]$, it follows that $F_u = 0$ if and only if $\operatorname{row}(F)\subseteq \operatorname{row}[q_1 \,\, \ldots \,\, q_{n_o}] = \operatorname{row}(\mathcal O)$, which is equivalent to statement 2. Therefore, if statement 2 holds, then $F_u = 0$ and statement 4 always holds for the unobservable subsystem $(C_u,A_u;F_u)$. Now we only have to prove that statement 2 implies statement 4 for the observable subsystem $(C_o,A_o;F_o)$. Since $(C_o,A_o)$ is full-state observable, it follows that
\begin{equation*}
    \rank
    \begin{bmatrix}
    \lambda_k I - J_k \\ C_k
    \end{bmatrix}
    = n_k, \quad \forall k,
\end{equation*}
\noindent
which implies statement 4.

\textit{(2)$\Leftarrow$(4)}. Consider the triple $(\bar C, J;\bar F)$. 
Let the submatrices $C_k$ and $F_k$ of $\bar C$ and $\bar F$, respectively, correspond to the Jordan submatrix $J_k$, for $k=1,\ldots,m$. The LHS of statement 4 can thus be expressed as
%
%
%
\def\arraystretch{1.1}
\begin{equation*}
    \rank
    \begin{bmatrix}
        \lambda I - J_1 & 0 & \ldots & 0 \\
        0 & \lambda I - J_2 & \ldots & 0 \\
        \vdots & \vdots & \ddots & \vdots \\
        0 & 0 & \ldots & \lambda I - J_m \\
        C_1 & C_2 & \ldots & C_m \\
        F_1 & F_2 & \ldots & F_m
    \end{bmatrix}
\end{equation*}
\noindent
and the LHS of statement 2 as
\begin{align*}
    \rank\begin{bmatrix} \mathcal O(\bar C,J) \\ \bar F \end{bmatrix} &= 
    \rank
    \begin{bmatrix}
        \mathcal O(C_1,J_1) & \ldots & \mathcal O(C_m,J_m) \\
        F_1 & \ldots & F_m
    \end{bmatrix}.
\end{align*}
\noindent
Due to the block diagonal structure of $J$, statement 4 holds if and only if
\begin{equation}
    \rank\begin{bmatrix}
        \lambda_k I - J_k \\ C_k \\ F_k
    \end{bmatrix}
    =
    \rank\begin{bmatrix}
        \lambda_k I - J_k \\ C_k
    \end{bmatrix}, \quad\forall k,
\label{eq.proof.pbhjordan}
\end{equation}
\noindent 
and statement 2 holds if and only if
\begin{equation}
    \rank\begin{bmatrix} \mathcal O(C_k,J_k) \\ F_k \end{bmatrix} = \rank(\mathcal O(C_k,J_k)), \quad\forall k,
\label{eq.proof.obsvjordan}
\end{equation}
\noindent
given that $\lambda_k$ is distinct for each submatrix $J_k$ and $\rank(C)=\rank(\bar C)=q$ without loss of generality.

Suppose statement 4 holds for the triple $(C,A;F)$ and, by similarity, condition \eqref{eq.proof.pbhjordan} also holds for all triples $(C_k,J_k;F_k)$. We now show that, under this assumption, condition \eqref{eq.proof.obsvjordan} holds for all $k$ and hence statement 2 also holds. If $F_k = 0$, conditions \eqref{eq.proof.pbhjordan} and \eqref{eq.proof.obsvjordan} always hold. Thus, in what follows we only consider the case $F_k\neq 0$.

If $J_k$ is a trivial Jordan block (with an eigenvalue of multiplicity $n_k=1$), condition \eqref{eq.proof.pbhjordan} is equivalent to
\begin{equation*}
    \rank\begin{bmatrix}
        0 \\ C_k \\ F_k
    \end{bmatrix}
    =
    \rank\begin{bmatrix}
        0 \\ C_k
    \end{bmatrix},
\end{equation*}
\noindent
which is satisfied if and only if $\row(F_k)\subseteq\row(C_k)$. The latter condition implies that condition \eqref{eq.proof.obsvjordan} is satisfied since $\mathcal O(C_k,J_k) = [(C_k)^\transp \,\, (\lambda_kC_k)^\transp \,\, \ldots \,\, (\lambda_k^{n-1}C_k)^\transp]^\transp$. Thus, Eq.~\eqref{eq.proof.pbhjordan} implies \eqref{eq.proof.obsvjordan} for trivial Jordan blocks.

Now, consider the nontrivial case of a Jordan submatrix $J_k$ corresponding to an eigenvalue of algebraic multiplicity $n_k>1$. For simplicity, and without loss of generality, consider the following Jordan submatrix $J_k$ and the corresponding submatrices $C_k$ and $F_k$:
\begin{equation*}
\begin{aligned}
\label{eq.proof.jordanblockexamp}
    J_k &= \operatorname{diag}(J_{k1},J_{k2},J_{k3}) 
    = \left[\begin{array}{ccc|cc|c}
        \lambda_k & 1 & 0 & 0 & 0 & 0 \\
        0 & \lambda_k & 1 & 0 & 0 & 0 \\
        0 & 0 & \lambda_k & 0 & 0 & 0 \\
        \hline
        0 & 0 & 0 & \lambda_k & 1 & 0 \\
        0 & 0 & 0 & 0 & \lambda_k & 0 \\
        \hline
        0 & 0 & 0 & 0 & 0 & \lambda_k
    \end{array}\right], \\
    C_k & = \left[\begin{array}{ccc|cc|c} \bm c_1 & \bm c_2 & \bm c_3 & \bm c_4 & \bm c_5 & \bm c_6 \end{array}\right], \\
    F_k & = \left[\begin{array}{ccc|cc|c} \bm f_1 & \bm f_2 & \bm f_3 & \bm f_4 & \bm f_5 & \bm f_6 \end{array}\right]. \\
\end{aligned}
\end{equation*}
\noindent
Condition \eqref{eq.proof.pbhjordan} 
can thus be written in terms of the rank of the following matrices
\def\arraystretch{1}
\begin{align*} \label{eq.proof.phbformjordan}
    \begin{bmatrix}
        \lambda_k I - J_k \\ C_k \\ F_k
    \end{bmatrix}
    &=
    \begin{bmatrix}\begin{smallmatrix}
        0 & 1 & 0 & 0 & 0 & 0 \\
        0 & 0 & 1 & 0 & 0 & 0 \\
        0 & 0 & 0 & 0 & 0 & 0 \\
        0 & 0 & 0 & 0 & 1 & 0 \\
        0 & 0 & 0 & 0 & 0 & 0 \\
        0 & 0 & 0 & 0 & 0 & 0 \\
        \bm c_1 & \bm c_2 & \bm c_3 & \bm c_4 & \bm c_5 & \bm c_6 \\
        \bm f_1 & \bm f_2 & \bm f_3 & \bm f_4 & \bm f_5 & \bm f_6
        \end{smallmatrix}
    \end{bmatrix}.
\end{align*}

\noindent
Since condition \eqref{eq.proof.pbhjordan} is satisfied, it follows that $[\bm f_1 \,\, \bm f_4 \,\, \bm f_6] = L_1 [\bm c_1 \,\, \bm c_4 \,\, \bm c_6]$ for some matrix $L_1\in\R^{r\times q}$. Given that the lead columns $\{\bm f_1,\bm f_4,\bm f_6\}$ of the Jordan blocks $\{J_{k1},J_{k2},J_{k3}\}$ are LI by assumption, then $\{\bm c_1,\bm c_4,\bm c_6\}$ are also LI. It thus follows from Lemma \ref{lemma.jordanleadingcols} that the pair $(C_k,J_k)$ is observable and hence $\rank(\mathcal O(C_k,J_k)) = n_k$. Thus, condition \eqref{eq.proof.obsvjordan} holds for nontrivial Jordan blocks. \QEDS
\end{pf}

\begin{cor} \label{cor.Adiagobsv}
    Suppose $A$ is diagonalizable. The triple $(C,A;F)$ is functionally observable if and only if
    \begin{equation}
    \label{eq.pbhrankfunctobsv}
        \rank \begin{bmatrix} \lambda I_n - A \\ C \\ F \end{bmatrix}
        =
        \rank \begin{bmatrix} \lambda I_n - A \\ C \end{bmatrix}, \,\, \forall\lambda\in\C.
    \end{equation}
\end{cor}

\begin{pf}
    \def\arraystretch{1}
    Consider the equivalent system $(\bar C,J;\bar F)$, where $J=PAP^{-1}$ is in Jordan form. If $A$ is diagonalizable, conditions \eqref{eq.proof.pbhjordan} and \eqref{eq.proof.obsvjordan} are respectively given by
    \begin{align}
        \rank\begin{bmatrix}
        C_k \\ F_k
    \end{bmatrix}
    &=
    \rank(C_k),
    \label{eq.proof.pbhjordan.diag}
    \\
    \rank\begin{bmatrix} C_k \\ \lambda_k C_k \\ \vdots \\ \lambda_k^{n_k-1} C_k \\ F_k \end{bmatrix} &= \rank\begin{bmatrix} C_k \\ \lambda_k C_k \\ \vdots \\ \lambda_k^{n_k-1} C_k \end{bmatrix}, 
    \label{eq.proof.obsvjordan.diag}
    \end{align}
    \noindent
    for $k=1,\ldots, m$, where $m$ is the number of distinct eigenvalues. Since conditions \eqref{eq.proof.pbhjordan.diag} and \eqref{eq.proof.obsvjordan.diag} are equivalent for all $k$, it follows that condition \eqref{eq.pbhrankfunctobsv} is equivalent to statement 2 in Theorem \ref{thm.functobsv}. \QEDS 
\end{pf}

The definition and conditions for functional observability have been established throughout several papers in the literature \citep{Darouach2000,Moreno2001,Jennings2011,Rotella2016,Darouach2023functdetect,Zhang2023onfunctobsv}. In the following remark, we clarify some of the recent development in functional observability and credit the individual contributions.

\textbf{Historical Remark 1}
    The notion of functional observability was originally proposed by \citet{Fernando2010} as a condition for the design of (functional) observers capable of asymptotically estimating the target state $\bm z(t)$ for arbitrary initial conditions. However, as we show in the equivalence between statements 1 and 2 in Theorem \ref{thm.functobsv}, functional observability is also more related to the unique reconstruction of $\bm z(0)$ from signals $\bm y(t)$ and $\bm u(t)$, without necessarily requiring the design of observers. This provides a result that is analogous to the case of full-state observability and the reconstruction of the full state $\bm x(t)$ using solely the observability Gramian. 

    \vspace{-0.1cm}
    Conditions analogous to Kalman's rank test and the PBH test for full-state observability were generalized to functional observability in  \citep{Jennings2011} and \citep{Moreno2001,Jennings2011}, respectively, leading to statements 2 and 4 in Theorem \ref{thm.functobsv}. The equivalence between statement 2 and 3 was later proven by \citet{Rotella2016}. Recently, \citet{Zhang2023onfunctobsv} highlighted that statement 4 is only necessary for functional observability and that sufficiency is only attained for special cases, such as for diagonalizable matrices $A$. Here, we show that the equivalence between statements 2 and 4 is determined by the Jordan form of a triple $(C,A;F)$, which allows us to relax the requirement of diagonalizability. This is illustrated in Example \ref{examp.cond4isonlynecessary} as follows. 

\begin{exmp}
\label{examp.cond4isonlynecessary}
    Consider the triple $(A,B;F)$ defined by
    \begin{equation} \label{eq.example2}
            A = \begin{bmatrix}
                0 & 1 & 0 \\
                0 & 0 & 1 \\
                0 & 0 & 0 
            \end{bmatrix},
            \quad
            C = \begin{bmatrix}
                0 & 0 & 1
            \end{bmatrix},
            \quad
            F = \begin{bmatrix}
                0 & 1 & 0
            \end{bmatrix},
    \end{equation}
    \noindent where $A$ is a $3\times 3$ Jordan block. By Definition \ref{def.jordanleadingcols}, the matrices $A$ and $F$ form a pair associated with the sole eigenvalue $\lambda_1 = 0$ of algebraic multiplicity $n_1=3$, where $F_{13}$ is the lead column of $F$ (which is a 1-dimensional vector in this case). The observability matrix is given by
    \begin{equation} \label{eq.examp.obsvmatrix}
        \mathcal O(C,A) = \begin{bmatrix}
                0 & 0 & 1 \\
                0 & 0 & 0 \\
                0 & 0 & 0 
            \end{bmatrix}.
    \end{equation}
    \noindent
    It follows from statement 2 in Theorem \ref{thm.functobsv} that the system is not functionally observable since $\rank([\mathcal O^\transp \,\, F^\transp]^\transp) = 2 \neq \rank(\mathcal O)$. On the other hand, statement 4 is satisfied since
    \begin{equation*}
        \rank \begin{bmatrix} A \\ C \\ F \end{bmatrix}
        =
        \rank \begin{bmatrix} A \\ C \end{bmatrix}
        =
        2
    \end{equation*}
    \noindent
    for $\lambda = 0$. This shows that statement 4 alone is not sufficient for functional observability. It is also clear that $F$ violates the assumption in Theorem \ref{thm.functobsv}: the lead column $F_{13} = 0$ and hence it does \textit{not} comprise a set of LI vectors. If we consider instead the functional matrix $F'=[1 \,\, 0 \,\, 0]$ (where $F'_{13}=1$ is the leading column), it can be verified that conditions 2 and 3 in Theorem \ref{thm.outputctrb} are not satisfied, which demonstrates the equivalence between these statements.
    \QEDA
\end{exmp}

The classical notion of full-state observability proposed by \citet{Kalman1959}, as well as the conventional PBH test \citep{Hautus1969}, follow as special cases of Theorem~\ref{thm.functobsv}.

\begin{cor}{\normalfont\textbf{(Full-state observability)}} Suppose $F = I_n$. The following statements are equivalent:
    \begin{enumerate}
        \item[1')] the triple $(C,A;I_n)$ is functionally observable;

        \item[2')] the pair $(C,A)$ is full-state observable;

        \item[3')] $\rank(\mathcal O) = n$;
        
        \item[4')] $\rank \begin{bmatrix} \lambda I_n - A \\ C \end{bmatrix} = n$, $\forall\lambda\in\C$.
    \end{enumerate}
\label{cor.fullobsv}
\end{cor}

\begin{pf}
    The equivalence between statements 1' and 2' follows from the fact that Definition~\ref{def.functobsv} reduces to the classical definition of full-state observability \citep[Definition 6.O1]{Chi-TsongChen1999} for $F=I_n$. Statements 2 and 4 in Theorem~\ref{thm.functobsv} reduce respectively to statements 3' and 4' since $\rank(F)=n$. For statement 4, note that $\rank(\bar F)=\rank(FP^{-1}) = \rank(F) = n$ and hence all columns of $\bar F$ are LI. \QEDS
\end{pf}

\rev{Note that the (generalized) PBH test can be directly applied to the system matrices $(C,A;F)$ without requiring their transformation to the Jordan form. This is the case provided that the assumption on the linear independence of the lead columns of $F_k$ is valid. 
This test establishes a condition based on the eigenstructure of the dynamical system that is  useful for the theoretical analysis of \textit{structured systems} \citep{ramos2022structural}. 
For example, \citet{Lin1974a} leveraged the PBH test to determine graph-theoretical conditions for \textit{structural} full-state controllability and, by duality, observability. By employing the generalized PBH test, \citet{Montanari2022} extended the structural observability conditions to the general case of structural functional observability. To this end, both papers considered classes of systems in which the PBH test is always valid. Specifically, \citet{Lin1974a} examined the special case $F=I_n$ (Corollary \ref{cor.fullobsv}), whereas \citet{Montanari2022} considered network systems with self-loops on all target variables (i.e, $A_{ii}\neq 0$ if $F_{ji}\neq 0$ for some $j=1,\ldots,r$), ensuring that all lead columns of $F_k$ are LI for a generic choice of parameters.
Later, \citet{Zhang2023onfunctobsv} proposed conditions for structural functional observability of general systems.
These grah-theoretical conditions for observability analysis have been pivotal in the development of cost-effective algorithms for the optimal placement of sensors in large-scale systems \citep{Liu2011,Montanari2022,eSousa2022,Montanari2023target,Zhang2023onfunctobsv,Zhang2023b}.
}

\section{\rev{PBH test for output controllability}}
\label{sec.outputctrb}

The output controllability of the dynamical system \eqref{eq.dynsys}--\eqref{eq.target} establishes the minimal conditions for the control of the target vector $\bm z(t)=F\bm x(t)$, without necessarily attaining full-state controllability. Since 1960, an algebraic rank condition for output controllability is available in the literature, providing a straightforward test based on matrix multiplications that generalizes Kalman's famous controllability test \citep{Bertram1960}. However, a test for output controllability based on the eigenspace of matrix $A$ (which is akin to the PBH test for full-state controllability) was missing, having only recently been established by \citet{Schonlein2023outputctrb}. As we show below in Example \ref{examp.cond3isonlynecessary}, the proposed PBH test for output controllability does not hold for \textit{all} dynamical systems, highlighting some gap in the derivation of this condition.
Given the recently established duality between functional observability and output controllability \citep{Montanari2023duality}, one might expect that the PBH condition for output controllability may be valid only under some assumption on the eigenstructure of the system (as we showed for the case of functional observability). Such derivation is presented next.

\begin{defn}
    The system \eqref{eq.dynsys}--\eqref{eq.target}, or the triple $(A,B;F)$, is output controllable if, for any initial state $\bm x(0)$ and any final target state $\bm z(t_1)$, there exists an input $\bm u(t)$ that steers $\bm z(0)=F\bm x(0)$ to $\bm z(t_1) = F\bm x(t_1)$ in finite time $t_1>0$.
\label{def.outputctrb}
\end{defn}

\begin{thm}
\label{thm.outputctrb}
{\normalfont\textbf{(Output controllability)}}
    Consider the triple $(A,B;F)$ and let $\mathcal E_i = \{\bm v\in\R^n \, : \, A^\transp\bm v = \lambda_i\bm v\}$ be the eigenspace of $A^\transp$ associated with $\lambda_i$, where $m$ is the number of distinct eigenvalues. The triple $(A,B;F)$ is output controllable if and only if one of the following conditions hold:
    \begin{enumerate}
        \item[1)] $\rank(F\mathcal C) = \rank(F)$;
        \smallskip

        \item[2)] $\rank \left(F\begin{bmatrix} \lambda I_n - A & B\end{bmatrix}\right) 
        =
        \rank (F)$, $\forall\lambda\in\C$, under the assumption that $(A,B)$ is full-state controllable;

        \item[3)] $\rank \left(F\begin{bmatrix} \lambda I_n - A & B\end{bmatrix}\right) 
        =
        \rank (F)$, $\forall\lambda\in\C$, and \begin{equation}
            \ker(\mathcal C^\transp)\cap\row(F)\cap\left\{ \bigcup_{i=1}^m \mathcal E_i\right\}\neq\emptyset,
        \label{eq.intersection}
        \end{equation}
        
        \noindent
        under the assumption that $(A,B)$ is not full-state controllable.
    \end{enumerate}
\end{thm}

\begin{pf}
%
%
\textit{Sufficiency and necessity of condition (1)}. Consider the input signal
\begin{equation}
\label{eq.proof.gramiancontrol}
    \bm u(t) = -B^\transp e^{A^\transp(t_1-t)} F^\transp W_F^{-1}(t_1) \bm (Fe^{At_1}\bm x(0) - \bm z(t_1)),
\end{equation}

\noindent
where $W_F(t_1) = F W_c(t_1)F^\transp$ is a projection of the controllability Gramian \rev{$W_c(t) = \int_0^{t} e^{A\tau}BB^\transp e^{A^\transp \tau}{\rm d}\tau$}. The response of the target $\bm z(t)$ at time $t=t_1$ is given by
\begin{equation}
\label{eq.proof.ztimeresponse}
    \begin{aligned}
    \bm z(t_1) &= 
    Fe^{At_1}\bm x(0) + F\int_0^{t_1} e^{A(t_1-t)}B\bm u(t) {\rm d}t \\
    &= Fe^{At_1}\bm x(0) - W_F(t_1) W_F^{-1}(t_1)\left(Fe^{At_1}\bm x(0) - \bm z(t_1)\right) \\
    &= \bm z(t_1).
    \end{aligned}
\end{equation}
\noindent
If \rev{$W_F(t_1)$} is invertible, then the system is output controllable according to Definition~\ref{def.outputctrb}. We now show the converse by contradiction. Suppose the system is output controllable but $W_F(t_1)$ is not invertible. Thus, there exists some nonzero $\bm v\in\R^{r}$ such that
\begin{equation*}
    \bm v^\transp FW_c(t_1)F^\transp\bm v = \int_0^{t_1} \norm{B^\transp e^{A^\transp (t_1-t)}F^\transp \bm v}^2 {\rm d}t = 0.
\end{equation*}
\noindent
\rev{By Definition \ref{def.outputctrb},} if $(A,B;F)$ is output controllable, there exists some $\bm u(t)$ that drives \rev{(in finite time $t_1$) any initial condition} $\bm z(0) = F\bm x(0) = Fe^{-At_1}\bar{\bm v}$ to \rev{the origin $\bm z(t_1) = 0$}, where $\bar{\bm v}$ is chosen such that $F\bar{\bm v}$ and $\bm v$ are linearly dependent. It then follows from Eq.~\eqref{eq.proof.ztimeresponse} that
%
\begin{equation*}
    \begin{aligned}
        0 &= F\bar{\bm v} + F\int_0^{t_1}e^{A(t_1-t)}B\bm u(t){\rm d}t,
        \\
        0 &= \bm v^\transp  F\bar{\bm v} + \bm v^\transp F\int_0^{t_1}e^{A(t_1-t)}B\bm u(t){\rm d}t,
        \\
        0 &= \bm v^\transp F\bar{\bm v},
    \end{aligned}
\end{equation*}
\noindent
which violates the assumption that $\bm v\neq 0$ since $\bm v$ and $F\bar{\bm v}$ are nonorthogonal.

To complete the proof, we show that the invertibility of $W_F(t_1)$ is equivalent to condition 1. Note that $\rev{\rank(FW_c(t_1)F^\transp) = r \Leftrightarrow} \row(FW_c)=\row(W_c)$, \rev{which is equivalent to condition 1, given that $\Im(W_c)=\Im(\mathcal C)$.}

\textit{(1)$\Leftrightarrow$(2)}. If $(A,B)$ is full-state controllable, then $\rank(\mathcal C)=n$ and $\rank[\lambda I_n-A \,\, B]=n$. Hence, condition 2 implies condition 1 trivially for any $F$. Contrariwise, if condition 1 does not hold, then $(A,B)$ is not full-state controllable, which is a contradiction to condition 2.

\textit{(1)$\Rightarrow$(3)}. We prove by contradiction. Suppose condition 1 is satisfied, but condition 3 is not. Thus, there exists an eigenvalue $\lambda_i$ and an $r$-dimensional vector $\bm q_i\neq 0$ such that $\bm q_i^\transp F[\lambda_iI_n - A \,\, B] = 0$. This implies that $\bm q_i^\transp FA = \lambda_i \bm q_i^\transp F$ and $\bm q_i^\transp FB = 0$. It thus follows that $\bm q_i^\transp FA^2 = \lambda_i\bm q_i^\transp FA = \lambda_i^2 \bm q_i^\transp F$ and, by induction, $\bm q_i^\transp FA^k = \lambda_i^k \bm q_i^\transp F$, for $k=1,2,\ldots$.
Multiplying $F\mathcal C$ on the left by $\bm q_i^\transp$ yields
\begin{equation*}
\begin{aligned}
    \bm q_i^\transp F\mathcal C &= 
    \bm q_i^\transp\begin{bmatrix}
        FB & FAB & \ldots & FA^{n-1} B
    \end{bmatrix} \\
    &= \begin{bmatrix}
        \bm q_i^\transp FB & \lambda_i \bm q_i^\transp FB & \ldots & \lambda_i^{n-1} \bm q_i^\transp FB
    \end{bmatrix} \\
    & = 0.
\end{aligned}
\end{equation*}
\noindent
This implies that $F\mathcal C$ is rank deficient, which contradicts the hypothesis that condition 1 is satisfied.

\textit{(1)$\Leftarrow$(3)}. 
We show that $\rank (\bar F\mathcal C) < r$ implies $\rank[\lambda  F - F A \,\,\, F B] < r$ for some $\lambda\in\C$. Under this condition, there exists some nonzero vector $\bm q\in\R^r$ that satisfies
\begin{equation}
\label{eq.proof.qFJBzero}
    \bm q^\transp F A^k B = 0, \quad {\rm for} \,\, k=0,1,2,\ldots .
\end{equation}
\noindent
It follows directly from Eq.~\eqref{eq.proof.qFJBzero} that $\bm q^\transp  F B = 0$. Therefore, in order to show that $\rank [F(\lambda I -  A) \,\,\, F B]<r$, it remains to prove that $\bm q^\transp F(\lambda I - A) = 0$, i.e., $F^\transp \bm q$ is an eigenvector of $A^\transp$.

The null space of $\mathcal C^\transp$ is defined by the subspace
\begin{equation*}
    \ker(\mathcal C^\transp) = \{ \bm w\in\R^n \, : \, \bm w^\transp A^kB = 0, \,\, \text{for} \,\, k=0,1,2,\ldots\}.
\end{equation*}
\noindent
Note that $\ker(\mathcal C^\transp)$ is non-empty since by assumption statement 2 does not hold and equivalently $(A,B)$ is not full-state controllable.
It thus follows from Eq.~\eqref{eq.proof.qFJBzero} that there exists some $\bm q\in\R^r$ such that $F^\transp\bm q\in\ker(\mathcal C^\transp)$. Therefore, 
\begin{equation}
\label{eq.proof.rowFkerC}
    \col(F^\transp) = \row(F)\cap\ker(\mathcal C^\transp) \neq \emptyset.
\end{equation}
\noindent
Condition \eqref{eq.intersection} is thus necessary and sufficient for the existence of some $F^\transp\bm q$ that lies simultaneously in $\ker(\mathcal C^\transp)$ and in some eigenspace of $A^\transp$, which implies that $F^\transp\bm q$ is an eigenvector. \QEDS
\end{pf}

\begin{figure}
    \centering
    \includegraphics[width=0.9\columnwidth]{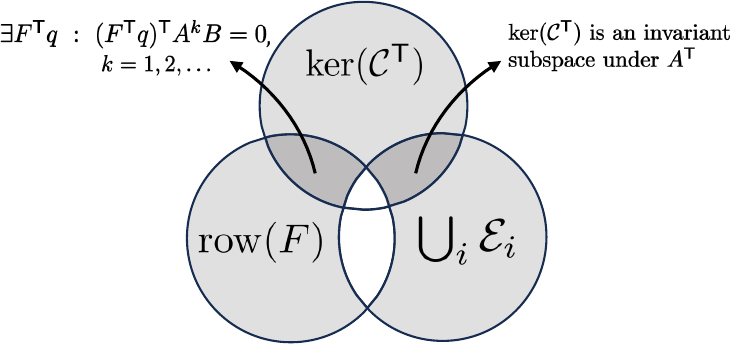}
    \caption{\label{fig.setintersection} Intersection between the row space of $F$, the null space of $\mathcal C^\transp$, and the union of all eigenspaces of $A^\transp$. The dark gray intersections are non-empty for all triples $(A,B;F)$ that are not full-state controllable. The white intersections can be empty for general triples $(A,B;F)$ (including full-state controllable and uncontrollable systems).}
\end{figure}

Figure \ref{fig.setintersection} illustrates the set intersection in condition \eqref{eq.intersection}. We further explore  condition \eqref{eq.intersection} to determine special classes of matrices $A$, $B$, and $F$ for which it always hold.

\begin{lem}
\label{lem.kercontainseigenvector}
    A non-empty subspace $\ker(\mathcal C^\transp)$ contains an eigenvector of $A^\transp$.
\end{lem}

\begin{pf}
    As proven in \citep{Schonlein2023outputctrb}, the subspace $\ker(\mathcal C^\transp)$ is invariant under the linear operator $A^\transp$. To see this, note that $A^\transp\bm w\in\ker(\mathcal C^\transp)$ for any $\bm w\in\ker(\mathcal C^\transp)$ because $(A^\transp\bm w)^\transp A^k B = \bm w^\transp A^{k+1} B = 0$. 
    Since $\ker(\mathcal C^\transp)$ is a nonempty invariant subspace under $A^\transp$, it follows that $\ker(\mathcal C^\transp)$ contains an eigenvector of $A^\transp$. \QEDS
\end{pf}

\begin{cor}
    \label{cor.rowFkerC}
    Let $\row(F)\supseteq\ker(\mathcal C^\transp)$. The triple $(A,B;F)$ is output controllable if and only if 
    \begin{equation}
        \label{eq.generalizedPBHcontrol}
        \rank(F[\lambda I_n - A \,\, B]) = \rank(F),  \,\, \forall\lambda\in\C.
    \end{equation}
\end{cor}

\begin{pf}
    For full-state controllable systems, the proof follows directly from condition 2 in Theorem \ref{thm.outputctrb}. We now prove this corollary for systems that are not full-state controllable [i.e., $\ker(\mathcal C^\transp)\neq\emptyset$].
    It follows from Lemma \ref{lem.kercontainseigenvector} that $\ker(\mathcal C^\transp)\cap\left\{ \bigcup_{i=1}^m \mathcal E_i\right\}\neq\emptyset$. Therefore, if $\row(F)\supseteq\ker(\mathcal C^\transp)$, then condition \eqref{eq.intersection} holds. \QEDS
\end{pf}

\begin{cor}
\label{cor.Adiagctrb}
    Let $A$ be diagonalizable. The triple $(A,B;F)$ is output controllable if and only if condition \eqref{eq.generalizedPBHcontrol} is satisfied.
\end{cor}

\begin{pf}
    If $A$ is diagonalizable, then $\bigcup_{i=1}^m \mathcal E_i = \R^n$. Since Eq. \eqref{eq.proof.rowFkerC} holds and $\row(F)\cap\ker(\mathcal C^\transp)\subseteq\R^n$, it follows that condition \eqref{eq.intersection} is always satisfied for diagonalizable matrices. \QEDS
\end{pf}

\begin{cor}{\normalfont\textbf{(Full-state controllability)}} Suppose $F = I_n$. The following statements are equivalent:
    \begin{enumerate}
        \item[1')] the triple $(A,B;I_n)$ is output controllable;

        \item[2')] the pair $(A,B)$ is full-state controllable;

        \item[3')] $\rank(\mathcal C) = n$;
        
        \item[4')] $\rank \begin{bmatrix} \lambda I_n - A & B \end{bmatrix} = n$, $\forall\lambda\in\C$.
    \end{enumerate}
\end{cor}

\begin{pf}
    The equivalence between statements 1' and 2' follows from the fact that Definition~\ref{def.outputctrb} reduces to the classical definition of full-state controllability \citep[Definition 6.1]{Chi-TsongChen1999} for $F=I_n$. Condition 1 in Theorem \ref{thm.outputctrb} directly reduces to statement 3' for $F=I_n$. Finally, if $F=I_n$, then $\row(F)=\R^n$. This implies that $\row(F)\supseteq\ker(\mathcal C^\transp)$ and hence the proof follows from Corollary \ref{cor.rowFkerC}. \QEDS
\end{pf}

\textbf{Historical Remark 2}
    The sufficiency and necessity of condition 1 was originally presented by \citet{Bertram1960} and is discussed in several textbooks (e.g., \citep[Section 9.6]{OgataBook}). Here, we presented an alternative derivation of this condition based on the controllability Gramian $W_c$ and its projection on the subspace spanned by the functional matrix $F$. The necessity of the generalized PBH rank condition \eqref{eq.generalizedPBHcontrol} was recently established by \citet{Schonlein2023outputctrb}. However, that paper inadvertently concludes that condition \eqref{eq.generalizedPBHcontrol} \textit{alone} is sufficient for \textit{any} triple $(A,B;F)$. Here, we revise the proof presented in \citep{Schonlein2023outputctrb} by showing that the generalized PBH condition \eqref{eq.generalizedPBHcontrol} requires the additional condition \eqref{eq.intersection} based on the eigenspace of $A$, the row space of $F$, and the null space of $\mathcal C^\transp$.
    The insufficiency of condition \eqref{eq.generalizedPBHcontrol} is further demonstrated in the following example.

\begin{exmp} 
\label{examp.cond3isonlynecessary}
    Consider the triple $(A,B;F)$ defined by
    \begin{equation*}
            A = \begin{bmatrix}
                0 & 1 & 0 \\
                0 & 0 & 1 \\
                0 & 0 & 0 
            \end{bmatrix},
            \quad
            B = \begin{bmatrix}
                1 \\ 0 \\ 0
            \end{bmatrix},
            \quad
            F = \begin{bmatrix}
                0 & 1 & 0
            \end{bmatrix}.
    \end{equation*}

    \noindent 
    %
    It follows from condition 1 in Theorem \ref{thm.outputctrb} that the system is not output controllable since $F\mathcal C = 0$. On the other hand, condition \eqref{eq.generalizedPBHcontrol} is satisfied since $[\lambda F-FA \,\, FB] = [0 \,\, 0 \,\, 1 \,\, 0]$ for $\lambda = 0$, yielding $\rank[0 \,\, 0 \,\, 1 \,\, 0] = \rank(F) = 1$. This shows that the generalized PBH condition \eqref{eq.generalizedPBHcontrol} alone is not sufficient for output controllability. It is also clear that $F$ violates condition \eqref{eq.intersection} in Theorem \ref{thm.outputctrb} given that $\row(F)$ does not intersect the eigenspace $\mathcal E_1 = \operatorname{span}\{[0 \, \, 0 \,\, 1]^\transp\}$ of $A^\transp$. If we consider instead the functional matrix $F'=[0 \,\, 0 \,\, 1]$, we can thus verify that conditions 1 and 3 in Theorem \ref{thm.outputctrb} are both satisfied, which demonstrates their equivalence.
    \QEDA
\end{exmp}

\section{Duality}
\label{sec.duality}

We present the duality between functional observability and output controllability, highlighting the relationship between the conditions established in Theorems \ref{thm.functobsv}~and~\ref{thm.outputctrb}. Note that the target vector \eqref{eq.target} has different interpretations depending on the considered problem. In the output controllability problem, $\bm z$ defines the vector sought to be controlled/steered, whereas, in the functional observability problem, $\bm z$ defines the vector sought to be estimated/reconstructed. We use $(A,B;F)$ to analyze the output controllability of a system described by the system matrix $A$, input matrix $B$, and functional matrix $F$, and 
$(C,A;F)$ to analyze the functional observability of a system described by the output matrix $C$, system matrix $A$, and functional matrix $F$.

The weak and strong duality between functional observability of the \textit{primal} system $(C,A;F)$ and the output controllability of the \textit{dual} (transposed) system $(A^\transp,C^\transp;F)$ is presented. 
%
\begin{prop}
\label{thm.weak}
{\normalfont\textbf{(Weak duality)}}
    If the primal system $(C,A;F)$ is functionally observable, then the dual system $(A^\transp,C^\transp;F)$ is output controllable.
\end{prop}

\begin{pf}
%
Recall that $\mathcal C(A^\transp,C^\transp) = \mathcal O(A,C)^\transp$. Therefore, statement 2 in Theorem \ref{thm.outputctrb} can be expressed as $\rank(F\mathcal O^\transp)=\rank(F)$. Since statement 2 in Theorem \ref{thm.functobsv} holds, it follows that $\row(F)\subseteq\row(\mathcal O)$ and, hence, $\rank(F\mathcal O^\transp) = \rank(F)$. \QEDS
%
%
\end{pf}

\begin{prop}
\label{thm.strong}
{\normalfont\textbf{(Strong duality)}}
    The primal system $(C,A;F)$ is functionally observable if and only if the dual system $(A^\transp,C^\transp;F)$ is output controllable and $F\Im(W)\perp F\ker(W)$, where $W = \int_0^{t_1} e^{A^\transp t}C^\transp Ce^{At}{\rm d}t$ is the observability Gramian of $(C,A)$ and, equivalently, the controllability Gramian of $(A^\transp,C^\transp)$.
\end{prop}

\begin{pf}
Since $\Im(W)=\Im(\mathcal C(A^\transp,C^\transp))=\Im(\mathcal O^\transp)$, it follows that if statement 2 in Theorem \ref{thm.outputctrb} holds, then $F\Im(W) = F\Im(\mathcal O^\transp) = \R^r$. This condition together with $F\Im(\mathcal O^\transp)\perp F\ker(\mathcal O^\transp)$ are equivalent to $\row(F)\subseteq\col(\mathcal O^\transp)=\row(\mathcal O)$, which in turn is equivalent to statement 2 in Theorem \ref{thm.functobsv}. \QEDS
\end{pf}

The generalized PBH tests for functional observability and output controllability are also related by the weak and strong duality principles. To see this, define 
    $\Psi := 
    \begin{bmatrix}
        \lambda I_n - A
        \\
        C
    \end{bmatrix}
    \,\,
    \text{and}
    \,\,
    \Psi^\transp =
    \begin{bmatrix}
        \lambda I_n - A^\transp & C^\transp
    \end{bmatrix}.$~For a pair of systems $(C,A;F)$ and $(A^\transp,C^\transp;F)$, we have that $\rank([\Psi^\transp \,\, F^\transp]^\transp) = \rank(\Psi)$ implies $\rank(F\Psi^\transp) = \rank(F)$, and the converse holds if and only if $F\operatorname{row}(\Psi) \perp F\ker(\Psi)$.

\textbf{Historical Remark 3}
    In this paper, we establish the duality between functional observability and output controllability based on the relation between their corresponding rank conditions. However, this duality can also be more generally defined from a geometric perspective: \citet[Section III]{Montanari2023duality} established the relationship between the set of functionally observable states $\mathcal O_F$ of a triple $(C,A;F)$ and the set of output controllable states $\mathcal C_F$ of the dual triple $(A^\transp,C^\transp;F)$. Specifically, the weak duality shows that $\mathcal O_F\subseteq \mathcal C_F$ holds for all pairs of systems, whereas the strong duality shows that $\mathcal O_F=\mathcal C_F$ holds if and only if $F\Im(W)\perp F\ker(W)$.

\rev{The following example illustrates the duality between these generalized properties, demonstrating that they are not always equivalent for a pair of transposed systems as in the special case of full-state controllability/observabilility.}

\begin{exmp}
    \rev{Consider matrices $A$ and $C$ in Eq.~\eqref{eq.example2} and $F = [1 \,\, 1 \,\, 1]$. Recall that $\mathcal O(C,A) = \mathcal C(A^\transp,C^\transp)^\transp$, where $\mathcal O(C,A)$ is given by Eq.~\eqref{eq.examp.obsvmatrix}. Following Theorems~\ref{thm.functobsv} and \ref{thm.outputctrb}, the system $(C,A;F)$ is not functionally observable, whereas the dual system $(A^\transp,C^\transp;F)$ is output controllable. Indeed, strong duality does not hold since $F\operatorname{Im}(\Psi) = F\operatorname{ker}(\Psi)= \operatorname{span}\{1\}$, which does not satisfy the orthogonality condition in Proposition \ref{thm.strong}. In contrast, for $F' = [1 \, 0 \, 0]$, strong duality holds given that $F'\operatorname{Im}(\Psi) = \operatorname{span}\{1\}$ and $F'\ker(\Psi) = \{0\}$, and hence $(C,A;F')$ is functionally observable and $(A^\transp,C^\transp;F')$ is output controllable.}
\end{exmp}

As we show next, the strong duality between functional observability and output controllability always hold for systems described by matrices $A$, $C$, and $F$ with special structures.

\begin{figure*}
    \centering
    \includegraphics[angle=0,origin=c,width=0.92\linewidth]{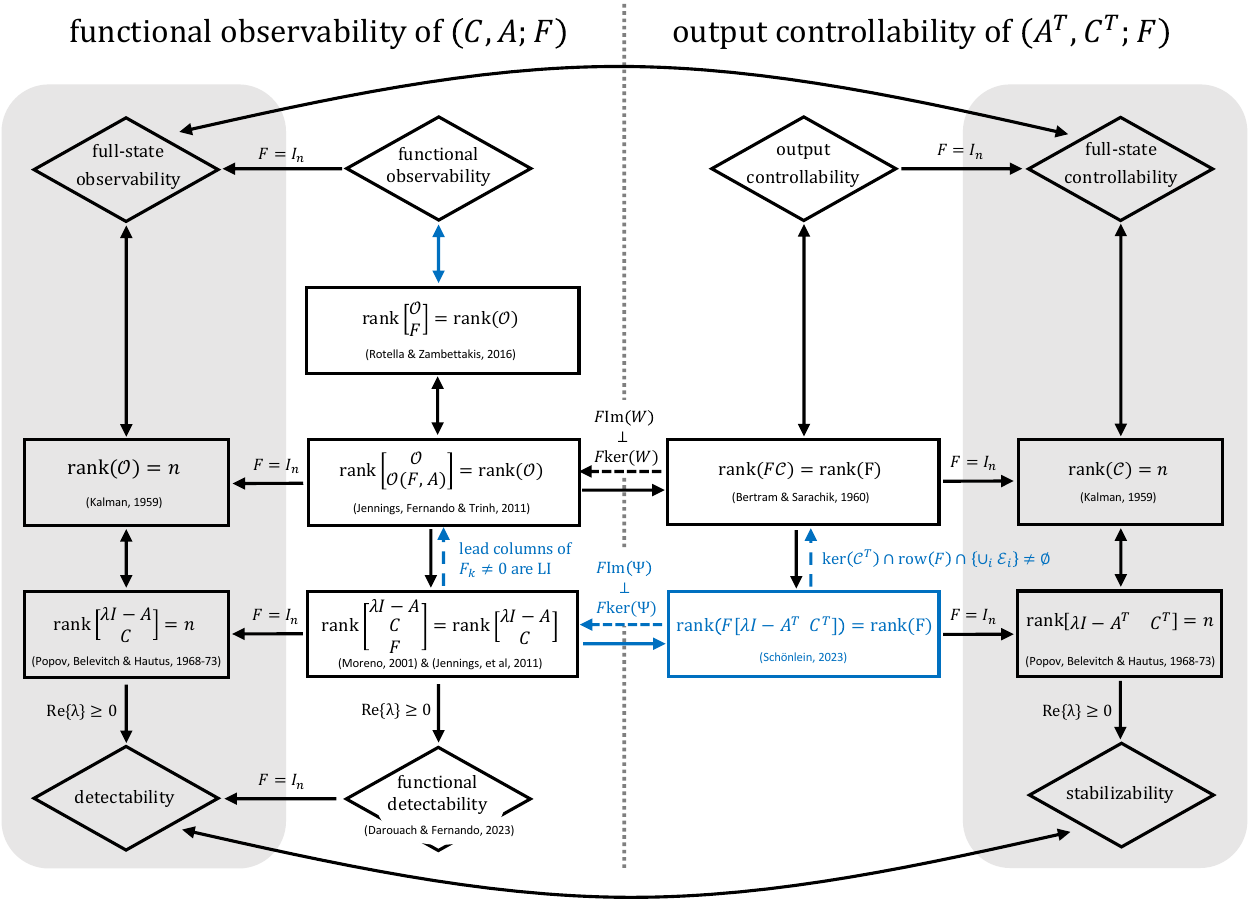}
    \caption{Relation between functional observability and output controllability. The square and diamond boxes represent conditions and definitions, respectively, while arrows indicate the equivalence or implication relations. 
    The novel conditions and relationships presented in this paper are highlighted as blue boxes and arrows, respectively. The matrices $\mathcal O$ and $\mathcal C$ correspond to the observability matrix of the pair $(C,A)$ and the controllability matrix of the pair $(A^\transp,C^\transp)$, respectively.}
\label{fig.diagram}
\end{figure*}

\begin{cor}
\label{cor.Adiagdual}
    Suppose that:
    \begin{enumerate}
        \item $A = U^{\transp}\Lambda U$, where $\Lambda = \diag(\{\lambda_i\}_{i=1}^n)$, $\lambda_i$ is an eigenvalue of $A$, and $U$ is a unitary matrix;
        \item the columns of $CU^\transp$ are orthogonal; and
        \item $\row(F)\subseteq\{\bigcup_{i=1}^n \mathcal E_i\}$.
    \end{enumerate}
    The primal system $(C,A;F)$ is functionally observable if and only if the dual system $(A^\transp,C^\transp;F)$ is output controllable.
\end{cor}

\begin{pf}
    Consider the triple $(C,A;F)$ and the equivalent system $(\bar C,\Lambda;\bar F)$, where $\bar C = CU^\transp$ and $\bar F = FU^{\transp}$. The observability Gramian of $(C,A)$ is thus given by
    \begin{equation*}
            W 
            =
            \int_0^{t_1} U^{\transp} e^{\Lambda t}U C^\transp C U^{\transp} e^{\Lambda t}U {\rm d}t \\
            = U^\transp\bar W U,
    \end{equation*}
    \noindent
    where $\bar W = \int_0^{t_1} e^{\Lambda t}\bar C^\transp \bar C e^{\Lambda t} {\rm d}t$. Since the columns of $\bar C$ are orthogonal by assumption, it follows that $\bar C^\transp \bar C$ is a diagonal matrix and hence $\bar W$ is also diagonal. It thus follows that $F\ker(W) = FU \ker(\bar W)$.
    By assumption, $\row(F)\subseteq \{\bigcup_{i=1}^n \mathcal E_i\} = \col(U)$. Therefore, $F\ker(W) = 0$, which trivially leads to $F\Im(W)\perp F\ker(W)$ and the strong duality (Proposition \ref{thm.strong}). \QEDS
\end{pf}

\begin{cor} {\normalfont\textbf{(Classical duality)}}
    Suppose $F = I_n$. The primal system $(C,A;F)$ is functionally observable if and only if the dual system $(A^\transp,C^\transp;F)$ is output controllable. 
\end{cor}

\begin{pf}
    The classical duality between full-state observability and controllability follows from Propositions \ref{thm.weak} and \ref{thm.strong} together with the fact that $\Im(W)\perp\ker(W)$. \QEDS
\end{pf}

\section{Conclusion}
\label{sec.conc}


The literature on functional observability and output controllability is vast and follows naturally from the 1960s work of Kalman on full-state observability and controllability. By examining the image space of the observability/controllability matrices (or their Gramians), it is clear that functional observability and output controllability directly generalize the corresponding notions of full-state observability and controllability. The conditions based on the observability/controllability matrices involve several matrix multiplications, which do not scale well with the system dimension. In contrast, the rank-based PBH conditions are directly tested on the system matrices (only requiring computation of their eigenvalues) and hence are better numerically conditioned for high-dimensional systems, such as large complex networks.
As a consequence, many applications to large-scale networks (e.g., optimal sensor placement, structured systems, attack/fault detection) are based on the generalized PBH rank condition \citep{emami2015functional,Montanari2022,eSousa2022,Venkateswaran2021,Zhang2023b}. 

However, the relationship between these properties is more subtle for conditions based on the eigenspace of the system matrix $A$, as given by the PBH conditions. Thus, the validity of the generalized PBH tests for functional observability and output controllability requires special assumptions on the matrices $A$, $C$ (or $B$), and $F$. We show that such assumptions are related to the Jordan form of these matrices, which hold trivially in the special case of diagonalizable matrices. As an immediate implication of our results, the generalized notion of \textit{functional detectability} \citep{Jennings2011,Darouach2023functdetect,Zhang2023onfunctobsv}\textemdash determined by the PBH test \eqref{eq.pbhrankfunctobsv} for unstable poles $\{\lambda\in\C \,:\,\operatorname{Re}\{\lambda\}\geq 0\}$\textemdash can now be applied to the more general class of systems satisfying the Jordan form assumed in Theorem \ref{thm.functobsv}.
Figure \ref{fig.diagram} summarizes the relation between the conditions for functional observability and output controllability presented in Sections \ref{sec.functionalobsv}--\ref{sec.duality}.

Our contributions are of relevance to theoretical results and methods based on the PBH conditions for functional observability and output controllability, especially in the context of large-scale systems.  
Such applications must take into consideration the classes of systems for which the generalized PBH conditions are valid, as proven in Theorems \ref{thm.functobsv} and \ref{thm.outputctrb}. 
We believe that our results will be important to stimulate research on generalizations and applications of these properties to broader classes of systems \citep{Trinh2006,Montanari2022nonlinear,muntwiler2023nonlinear}.

\begin{ack}                 
The authors acknowledge support from the U.S. Army Research Office, Grant No. W911NF-19-1-0383 (A.N.M. and A.E.M.), and the National Natural Science Foundation of China, Grant No. GQQNKP001 (C.D.).
\end{ack}



\begin{small}

\end{small}

\end{document}